%% file: plane_perc.tex
\newcommand\webcite[1]{\texttt{\def~{\~{}}#1}}
\documentclass{article}
\usepackage{amsmath,amsfonts,amsthm,epsfig,color}
\newtheorem{theorem}{Theorem}
\newtheorem{lemma}[theorem]{Lemma}

\newtheorem{corollary}[theorem]{Corollary}

\def\Sectionpfmain{Section 8}% of Voronoi. Check when published!
\def\TheoremVRSW{Theorem 12}% of Voronoi. Check when published!
\def\SectionsecS{Section 5}%sec_K of Kesten. Check when published!
\def\Lemmalongweak{Lemma~11}% Of Kesten. Check when published!
\def\Lemmalongstrong{Lemma~9}% Of Kesten. Check when published!
\def\TheoremHarris{Theorem 8}% of Kesten.
\def\LemmaHalf{Lemma 3}% of Kesten.
\def\Lemmahalf{Corollary~5}% Of Kesten. Check when published!
\def\LemmaX{Lemma~6}% Of Kesten. Check when published!

\def\pp#1{\Pr^\pi_{#1}}
\def\ppw{\Pr_p^w}
\def\phw{\Pr_{1/2}^w}
\def\pqw{\Pr_{1-p}^w}
\def\thpw{\theta^w(p)}
\def\LP{{LV}}
\def\Lt{{L_\triangle}}
\def\epsilon{{\varepsilon}} % GLOBAL REPLACE ??
\def\star{{*}} % GLOBAL REPLACE ??

\def\Prt{{\widetilde\Pr}}
\def\Prtwo{{\widetilde\Pr}}

\newcommand\E{{\mathop{\mathbb E{}}\nolimits}}
\renewcommand\Pr{{\mathop{\mathbb P{}}\nolimits}}

\def\Ls{{L_\square}}

\def\Ldual{{L^\star}}
\def\Lsd{{L_{\raise.05pt\rlap{$\mkern-1.5mu\times$}\square}}}
\def\Lsdsub{{L_{\raise-.40pt\rlap{$\mkern-2.2mu\times$}\square}}}%BROKEN!!
\newcommand{\cA}{{\mathcal A}}
\newcommand{\cB}{{\mathcal B}}
\newcommand{\Pow}{{\mathcal P}}%mathop etc?
\newcommand{\Z}{{\mathbb Z}}
\newcommand{\hZs}{(\frac{1}{2}{\mathbb Z})^2}
\newcommand{\RR}{{\mathbb R}}

\begin{document}
\title{Sharp thresholds and percolation in the plane}
\date{October 4, 2005}

\author{B\'ela Bollob\'as\thanks{Department of Mathematical Sciences,
University of Memphis, Memphis TN 38152, USA}
\thanks{Trinity College, Cambridge CB2 1TQ, UK}
\thanks{Research supported in part by NSF grant ITR 0225610 and DARPA grant
F33615-01-C-1900}
\thanks{Research partially undertaken during a visit to the Forschungsinstitut f\"ur Mathematik,
ETH Z\"urich}
\and Oliver Riordan$^{\dag\S}$%
\thanks{Royal Society Research Fellow, Department of Pure Mathematics
and Mathematical Statistics, University of Cambridge, UK}}
\maketitle

\begin{abstract}
Recently, it was shown in~\cite{Voronoi} that the critical probability for random
Voronoi percolation in the plane is $1/2$. As a by-product of the method, a short
proof of the Harris-Kesten Theorem was given in~\cite{ourKesten}. The aim of this paper
is to show that the techniques used in these papers can be applied
to many other planar percolation models, both to obtain short proofs of known results, and to
prove new ones.
\end{abstract}

\section{Introduction}

In~\cite{ourKesten}, a short proof was given of the fundamental result
of Harris~\cite{Harris} and Kesten~\cite{Kesten1/2} that the critical probability
$p_H=p_H(\Z^2,\mathrm{bond})$ for bond percolation in the planar square lattice $\Z^2$ is equal to $1/2$,
where $p_H$ is the critical probability for the occurrence of percolation (see below),
and $\Z^2$ is the graph with vertex set $\Z^2$ in which vertices are adjacent if and only
if they are at Euclidean distance $1$.
The methods used in~\cite{ourKesten} were developed in~\cite{Voronoi} to prove
the new result that the critical probability for percolation in random plane Voronoi
tilings is also $1/2$. Here we show that the same methods easily give exponential decay
of the volume below the critical probability.
Furthermore, while the arguments in~\cite{ourKesten} are written specifically
for bond percolation in $\Z^2$, they can also be applied in many other planar
contexts. We illustrate this by
considering several examples. We start with two well-known ones,
site percolation in the square and triangular lattices.
Next, we consider a new bond percolation model in the square lattice,
where the states of the edges are not independent, showing that an analogue
of the Harris-Kesten result holds in this context.
Finally, we study random discrete Voronoi percolation in the plane.
It is very likely that the methods of~\cite{Voronoi} and~\cite{ourKesten} can be
applied to many other percolation models.

In the rest of this introduction we shall recall some of the fundamental concepts of
percolation theory. Then, in Section~\ref{sec_prelim}, we present the basic tools
we shall use to prove our results. In Section~\ref{bdecay} we show that the method
of~\cite{ourKesten} easily extends to prove an exponential decay result of Kesten~\cite{Kesten81}.
In Section~\ref{sec_other} we apply our method to give short proofs of well-known results
for site percolation in the square and triangular lattices. Finally, in Section~\ref{sec_non}
we consider two percolation models that do not correspond to (independent) site percolation
on any lattice, proving results we believe to be new.

A {\em bond percolation measure} on an infinite graph $G$ is a probability
measure on the space of assignments of a {\em state}, namely {\em open} or {\em closed},
to each edge $e$ of $G$ (with the usual
$\sigma$-field of measurable events).
Similarly, a {\em site percolation measure} on $G$
is a probability measure on assignments of states to vertices. Here, $G$ will
usually be a planar lattice; in particular, we consider the square lattice $\Z^2$
and the triangular lattice $\Lt$.

Given a lattice $L$, when discussing bond percolation on $L$ we consider
the measure $\Pr_p^{L,\mathrm{bond}}$
in which the states of the edges are independent, and
each edge is open with probability $p$. Similarly, when discussing
site percolation on $L$ we consider the measure $\Pr_p^{L,\mathrm{site}}$ in which the
vertices are open independently with probability $p$. When there is no danger of confusion,
we write $\Pr_p$ for either of these measures.

An {\em open cluster} is a maximal connected subgraph of $L$ all of whose edges (vertices)
are open. We write $C_v$ for the open
cluster containing a given vertex $v\in L$. Thus a vertex $w$ lies in $C_v$
if and only if $w$ can be reached from
$v$ by an open path, i.e., a path in $L$ all of whose edges (vertices) are open.
In the case of site percolation, if $v$ is closed then $C_v=\emptyset$.

Writing $|C_v|$ for the number of vertices of $C_v$, let
\[
 \theta(p)=\Pr_p\left(|C_0|=\infty\right),
\]
where $0=(0,0)$ is the origin. We shall always take $0$ to be a vertex of $L$.
By Kolmogorov's $0$-$1$ law, percolation occurs if and only if $\theta(p)>0$. More precisely,
if $\theta(p)>0$, then with probability $1$ there is an infinite open cluster
somewhere in $L$, while if $\theta(p)=0$, then with probability $1$ there is no such cluster.
As $\theta(p)$ is increasing in $p$, there is a critical probability $p_H$
such that $\theta(p)=0$ for $p<p_H$ and $\theta(p)>0$ for $p>p_H$.
This critical probability depends on the lattice
$L$ and type of percolation under consideration. To emphasize this dependence we may write
$p_H(L,\mathrm{bond})$ or $p_H(L,\mathrm{site})$.
Here, following Welsh
%~\cite{WelshSciProg}
(see~\cite{SW}), the $H$ is in honour of Hammersley;
Broadbent and Hammersley introduced the basic concepts of percolation in a 1957
paper~\cite{BH}, where they posed the problem of determining $p_H$ in a variety of contexts.
Hammersley~\cite{H2,H4,H5} proved general upper and lower bounds which imply, for example,
that $0.35<p_H(\Z^2,\mathrm{bond})<0.65$.

Writing $\E_p$ for the expectation corresponding to $\Pr_p$, let
\[
 \chi(p) = \E_p |C_0|
\]
be the expected size of the open cluster of the origin.
It is immediate that $\chi(p)$ is increasing in $p$, so there is another critical probability,
\[
 p_T = \inf\{p : \chi(p)=\infty\},
\]
with the $T$ in honour of Temperley.
As $\theta(p)>0$ implies $\chi(p)=\infty$, we have $p_T \le p_H$.

For many years it was believed that $p_T=p_H=1/2$ for bond percolation in $\Z^2$;
this conjecture seems not be have been made explicitly, but, supported by
various results and numerical evidence, this belief gradually arose.
In 1978, Russo~\cite{Russo} and Seymour and Welsh~\cite{SW} made significant progress.
In particular, they proved independently
that $p_T+p_H=1$.
It was only in 1980, twenty years after Harris' proof of the inequality $p_H\ge 1/2$,
that Kesten~\cite{Kesten1/2} proved that $p_T=p_H=1/2$.
Since then, Menshikov~\cite{Menshikov} (see also Menshikov, Molchanov and Sidorenko~\cite{MMS})
and Aizenman and Barsky~\cite{AB} (see also Grimmett~\cite{Grimmett})
have shown that $p_T=p_H$ in great generality, in particular, for site percolation
in any lattice graph; see Section \ref{ss_sth} for a formal definition.
Note that bond percolation in a lattice graph $L$ corresponds to site percolation in
the line graph of $L$, which can be realized as a lattice graph, so results for site percolation
in general lattices apply to bond percolation as well.

Below the critical probability, much stronger
results are known than $\chi(p)<\infty$. In particular, Kesten~\cite{Kesten81} showed in 1981
that for site percolation in a lattice, when $p<p_T$, the number $|C_0|$ of vertices
in $C_0$ decays exponentially.
(See also Aizenman and Newman~\cite{AN} and Grimmett~\cite{Grimmett}.)
In the light of the proofs that $p_T=p_H$ mentioned above,
Kesten's result implies that there is a single critical probability $p_H$, with percolation above $p_H$
and exponential decay of the size of the open cluster of the origin below $p_H$.
Here we shall show that, in various contexts, the method of~\cite{ourKesten}
easily gives exponential decay for $p<p_H$, implying that $p_T=p_H$.

\begin{figure}[htb]
 \[\epsfig{file=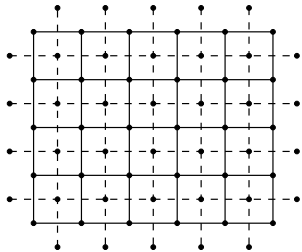,height=1.3in}\]
\caption{Portions of the lattice $L=\Z^2$ (solid lines) and the isomorphic dual lattice
$\Ldual$ (dashed lines).}
\label{fig_selfdual}
\end{figure}

An important property of bond percolation in $\Z^2$ is the `self-duality' of $\Z^2$. This property is key
to the results of Harris and Kesten. In the context of bond percolation, the appropriate notion of duality
is the standard one for plane graphs: the dual $G^\star$ of a graph $G$ drawn in the plane has a vertex for each
face of $G$, and an edge $e^\star$ for each edge $e$ of $G$. The edge $e^\star$ joins the two vertices of $G^\star$
corresponding to the faces of $G$ in whose boundary $e$ lies.
Taking $G=\Z^2$, there is a vertex $v$ of $G^\star$ for each square $[a,a+1]\times[b,b+1]$, $a,b\in \Z$,
which we may take to be the point $v=(a+1/2,b+1/2)$. It is easy to see that $G^\star$ is isomorphic to $G$; see Figure~\ref{fig_selfdual}.
This self-duality can be considered the `reason why' $p_H(\Z^2,\mathrm{bond})=1/2$, but 
this trivial observation, made soon after the question first arose,
is very far from giving a proof of the Harris-Kesten result.

\section{Preliminaries}\label{sec_prelim}

As in~\cite{ourKesten}, the proofs here will be mostly self-contained. The main
result we shall use is a sharp-threshold result of Friedgut and Kalai~\cite{FK},
a simple consequence of a result of Kahn, Kalai and Linial~\cite{KKL} concerning
the influences of coordinates in a product space. (See also~\cite{BKKKL}.)

Let $X$ be a fixed ground set with $N$ elements, and let $X_p$ be a random subset
of $X$ obtained by selecting each $x\in X$ independently with probability $p$.
For a family $\cA\subset \Pow(X)$ of subsets of $X$, let $\Pr_p^X(\cA)$ be the probability that $X_p\in\cA$.
In this context, $\cA$ is {\em increasing} if $A\in \cA$ and $A\subset B\subset X$
imply $B\in \cA$.
Also, $\cA$ is {\em symmetric} if there is a permutation group acting transitively
on $X$ which fixes $\cA$. In other words, $\cA$ is a union of orbits of the induced action
on $\Pow(X)$.
In our notation the result of Friedgut and Kalai~\cite{FK} we shall need is as follows.

\begin{theorem}\label{th_sharp}
There is an absolute constant $c_1$ such that if
$|X|=N$, the family
$\cA\subset \Pow(X)$ is symmetric and increasing, and
$\Pr_p^X(\cA)>\epsilon$, then $\Pr_q^X(\cA)>1-\epsilon$ whenever
$q-p \ge c_1\log(1/(2\epsilon))/\log N$.
\end{theorem}

We shall also make frequent use of Harris' Lemma.

\begin{lemma}\label{l_Harris}
If $\cA$, $\cB\subset \Pow(X)$ are increasing, then for any $p$ we have
\[
 \Pr_p^X(\cA\cap \cB)\ge \Pr_p^X(\cA) \Pr_p^X(\cB).
\]
\end{lemma}

Taking complements, the lemma also applies to two {\em decreasing} events, where a decreasing
event is the complement of an increasing one.
In other contexts, Lemma~\ref{l_Harris} is often known as Kleitman's Lemma~\cite{Kleitman}.
The present context is exactly that of Harris' original paper~\cite{Harris}: $X$ will be a set of edges
or vertices in the lattice (according to whether we are considering site or bond percolation), and
$X_p$ will be the subset of $X$ consisting of the open edges/vertices. Thus an event is increasing if it
is preserved by changing the states of one or more edges/vertices from closed to open, and Harris' Lemma
states that increasing events are positively correlated.

In addition to the results above, we shall need two observations concerning $k$-dependent percolation.
A bond percolation measure on a graph $G$ is {\em $k$-dependent} if,
for every pair $S$, $T$ of sets of edges of $G$ at graph
distance at least $k$, the states (being open or closed) of the edges in $S$ are independent of
the states of the edges in $T$. When $k=1$, the separation
condition is exactly that no edge of $S$ shares a vertex with an edge of $T$.
The definition of $k$-dependence for a site percolation measure on $G$ is exactly the same, except
that $S$ and $T$ run over all sets of {\em vertices} at graph distance at least $k$. Here we shall consider
dependent measures only on the lattice $\Z^2$.

These $k$-dependent measures arise very naturally in a variety of contexts (for example,
static renormalization arguments), and have been considered by several authors;
see Liggett, Schonmann and Stacey~\cite{LSS} and the references therein.
In~\cite{LSS}, a very general comparison result between $k$-dependent and product measures is proved:
working on any fixed countable graph $G$ of bounded degree (for example, $\Z^d$),
for any $p<1$ there is an $f(G,k,p)<1$ such that any $k$-dependent measure
in which each edge (vertex) is open with probability at least $f(G,k,p)$
dominates the product measure $\Pr_p$ in which edges (vertices)
are open independently with probability $p$.

In particular, provided the individual edge probabilities are high enough,
percolation occurs in $\Z^2$ under the assumption of $1$- (or $k$-)
dependence.

\begin{lemma}\label{l_1dep}
There is a $p_0<1$ such that in any $1$-dependent bond percolation measure on $\Z^2$ satisfying the
additional condition that each edge is open with probability at least $p_0$, the probability
that $|C_0|=\infty$ is positive.
\end{lemma}

In applications, the value of $p_0$ is frequently important. Currently, the best known bound is the
result of Balister, Bollob\'as and Walters~\cite{BBW} that one can take $p_0=0.8639$.
Here, the value of $p_0$ will be irrelevant: all we shall need is the
essentially trivial Lemma~\ref{l_1dep}.
For completeness, we give a very simple proof that one can take $p_0=0.995$.

Indeed, suppose that the open cluster $C_0$ containing the origin is finite.
Let $C_\infty$ be the (unique) infinite component of $\Z^2\setminus C_0$,
and let $B$ be the edge-boundary of $C_\infty$, i.e., the set of edges joining
$C_\infty$ to $\Z^2\setminus C_\infty$. Note that every edge in $B$ joins
$C_\infty$ to $C_0$, and is thus closed.
Passing to the lattice $L^\star$ dual to $L=\Z^2$ as defined above,
the edges of $L^\star$ corresponding to the edges of $L$ in $B$
form a simple cycle $S$ in $L^\star$ that surrounds the
origin.

\begin{figure}[htb]
 \[\epsfig{file=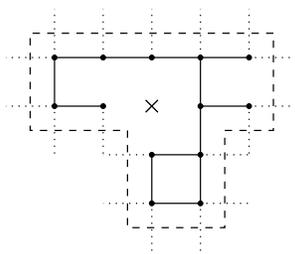,height=1.3in}\]
\caption{An open cluster $C$ in $L=\Z^2$ (dots and solid lines), the edge boundary $B$
of the infinite component $C_\infty$ of $L\setminus C$ (dotted lines), and the corresponding cycle $S$
in $L^\star$ (dashed lines). The point marked with a cross is in a finite
component of $L\setminus C$.}
\label{fig_dualcycle}
\end{figure}

Given the length $\ell\ge 4$ of $S$, there are crudely at most $\frac{\ell-2}{2}3^{\ell-2}$ possibilities
for $S$ (and hence $B$): $S$ must cross the $x$ axis at some $x$-coordinate between $\frac{1}{2}$ and $\frac{\ell-3}{2}$.
Walking round $S$, at each stage
there are at most three possibilities for the next edge, and at most one
choice that closes the cycle at the end.
Passing back to $L=\Z^2$,
the edges of $L$ may be partitioned into four complete matchings, one of which must contain a 
set $B'$ of at least $|B|/4=\ell/4$ edges of $B$.
Now the states of the edges in $B'$ are independent of each other, and each $e\in B'$
is closed with probability at most $1-p_0$.
Putting everything together, we see that the probability that $|C_0|$ is finite,
which is exactly the probability that some closed cycle in the dual surrounds the origin, is at most
\[
 \sum_{\ell\ge 4,\,\text{$\ell$ even}} \frac{\ell-2}{2} 3^{\ell-2} (1-p_0)^{\ell/4}.
\]
This is strictly less than $1$ if $p_0=0.995$.

Finally, a corresponding negative result is just as easy:
we repeat the statement and proof from~\cite{Voronoi}.
This time, it is easier to work with site percolation.
Recall that in the site percolation context, $C_0$, the open cluster of the origin,
is the set of vertices of $\Z^2$ joined to the origin by a path in $\Z^2$ every one of whose vertices is open.

\begin{lemma}\label{l_kdepneg}
Let $k$ be a fixed positive integer, and let $\Prtwo$ be a $k$-dependent site percolation measure on $\Z^2$
in which
every vertex $v\in \Z^2$ is open with probability at most $p$. There is a constant $p_1=p_1(k)>0$ such that
for every $p\le p_1$ there is a $c(p,k)>0$ for which
\[
 \Prtwo( |C_0| \ge n)\le \exp(-c(p,k)n)
\]
for all $n\ge 1$.
\end{lemma}
\begin{proof}
If $|C_0|\ge n$, then the subgraph of $\Z^2$ induced by the open vertices
contains a tree $T$ with $n$ vertices, one of which is the origin.
It is well known and easy to check that the number of such trees in $\Z^2$ grows exponentially, and is at most $(4e)^n$.
Fix any such tree $T$. Then there is a subset $S$ of at least $n/(2k^2-2k+1)$ vertices of $T$ such that any $a,b\in S$
are at graph distance at least $k$; indeed, one can find such a set by a greedy algorithm:
whenever a vertex $a$ is chosen, the number of other vertices it rules out is at most the number
of other vertices of $\Z^2$ within graph distance $k-1$ of $a$, namely $4\binom{k}{2}=2k^2-2k$.
The vertices of $S$ are open independently, so the probability that every vertex of $T$ is open is at
most $p^{|S|}$.
Hence,
\[
  \Prtwo( |C_0| \ge n)\le (4e)^n p^{n/(2k^2-2k+1)}.
\]
Provided $p$ is small enough that $r=4e p^{1/(2k^2-2k+1)}<1$, the conclusion follows,
taking $c(p,k)=-\log r$.
\end{proof}

\section{Bond percolation in $\Z^2$: exponential decay}\label{bdecay}

In this section we consider bond percolation in $\Z^2$,
writing $\Pr_p$ for the probability measure $\Pr_p^{\Z^2,\mathrm{bond}}$, in which 
each edge of $\Z^2$ is open with probability $p$, independently of all other edges.
In~\cite{ourKesten}, a short proof was given of the Harris-Kesten result that in this context $p_H=1/2$,
using Theorem~\ref{th_sharp} as the main ingredient.
In fact, the method also gives a simple proof that for $p<1/2$ there is
exponential decay of the `volume' $|C_0|$ of the open cluster containing the origin.
It follows that $\chi(p)$ is finite for $p<1/2$, and hence that $p_T=p_H=1/2$.
The result below was first proved by Kesten~\cite{Kesten81} in 1981.

\begin{theorem}\label{th_Lsdecay}
For every $p<1/2$, there is a constant $a=a(p)>0$ such that
$\Pr_p(|C_0| \ge n)\le \exp(-an)$ for all $n\ge 0$.
\end{theorem}

We shall deduce Theorem~\ref{th_Lsdecay} from \Lemmalongweak\ of~\cite{ourKesten},
reproduced below as Lemma~\ref{l_longhigh}.
Most of the work in~\cite{ourKesten} went into proving this lemma (or the stronger
form, \Lemmalongstrong\ in \cite{ourKesten});
the deduction of the Harris-Kesten
Theorem was then easy.
The lemma concerns `open crossings of rectangles':
we identify a rectangle $R=[x_0,x_1]\times [y_0,y_1]$, where $x_0<x_1$ and $y_0<y_1$ are integers,
with an induced subgraph of $\Z^2$. This subgraph includes all vertices and edges in the interior
and boundary of $R$. We write $H(R)$ for the event that there is a {\em horizontal open crossing} of $R$,
i.e., a path from the left side of $R$ to the right side consisting entirely of open edges
of $R$.
Similarly, we write $V(R)$ for the event that there is a vertical open crossing
of $R$.

\begin{lemma}\label{l_longhigh}
Let $p>1/2$ be fixed.
If $R_n$ is a $3n$ by $n$ rectangle in $\Z^2$, then
$\Pr_p(H(R_n))\to 1$ as $n\to\infty$.
\end{lemma}

\begin{proof}[Proof of Theorem~\ref{th_Lsdecay}]
Fix $p<1/2$, 
let $p_1>0$ be a constant for which Lemma~\ref{l_kdepneg} holds with $k=9$, and set $c=(1-p_1)^{1/4}$.

We shall apply Lemma~\ref{l_longhigh} to the lattice $\Ldual$ dual to $L=\Z^2$,
which is isomorphic to $\Z^2$.
Defining the state of a dual edge $e^\star$ to be the state of $e$, each edge of $\Ldual$ is closed with probability
$1-p>1/2$, independently of all other edges.
By Lemma~\ref{l_longhigh}, if $R$ is a $3m$ by $m$ rectangle
in $\Ldual$ then, provided we choose $m\ge 10$ large enough,
the probability that $R$ is crossed the long way by a path of closed dual
edges is at least $c$.

Set $s=m-1$, and let $S$ be an $s$ by $s$ square in $\Z^2$. Arrange four $3m$ by $m$
rectangles in the dual lattice to form an annulus $A$ as in Figure~\ref{fig_ann},
\begin{figure}[htb]
 \[\epsfig{file=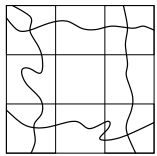,height=0.8in}\]
\caption{Four rectangles forming an annulus.}
\label{fig_ann}
\end{figure}
with the inside of the annulus surrounding $S$.
Using Lemma~\ref{l_Harris}, with probability at least $c^4=1-p_1$,
each of the four rectangles is crossed the long way by a path of closed
dual edges. If this happens, then there is a cycle of closed dual edges in $A$ which surrounds $S$.
(See Figure~\ref{fig_ann}.)
It follows that in this case,
in the original lattice, no vertex in $S$ is connected by an open path to a vertex outside $A$.

Returning to $\Z^2$, given an $s$ by $s$ square $S$ in $\Z^2$,
let $B(S)$ be the event that some vertex in $S$ is connected by an open path
to a vertex at $L_\infty$-distance $2s>m+1$ from $S$. We have shown that
$\Pr_p(B(S))\le p_1$.

Let us define a site percolation measure $\Prt$ on $\Z^2$ as follows:
each $v=(x,y)\in \Z^2$ is open if $B(S)$ holds
for the square $S_v=[sx,s(x+1)]\times [sy,s(y+1)]$. As $B(S_v)$ depends only on the states of edges within
$L_\infty$-distance $2s$ of $S_v$, the measure $\Prt$ is $9$-dependent. Furthermore, each
$v\in \Z^2$ is open with $\Prt$-probability at most $p_1$. Let $C_0$ be the open cluster of the origin
in our original bond percolation, and let $C_0'$ be the open cluster of the origin
in the site percolation we have just defined.
By Lemma~\ref{l_kdepneg} there is an $a>0$ such that
\[
 \Prt(|C_0'|\ge n)\le \exp(-an)
\]
for every $n$.

If $|C_0|>(6s+1)^2$, then every vertex $w$ of $C_0$ is joined
by an open path to some vertex at $L_\infty$-distance $3s$ from $w$. If $w\in S_v$, then it follows
that $B(S_v)$ holds. Thus, if $|C_0|>(6s+1)^2$, then $B(S_v)$ holds for every $v$ such that $S_v$ contains vertices
of $C_0$. The set of such $v$ forms an open cluster with respect to $\Prt$, and is thus a subset of $C_0'$.
Hence, as each
$S_v$ contains only $(s+1)^2$ vertices, for
$n\ge (6s+1)^2$ we have
\[
 \Pr_p(|C_0|\ge n) \le \Prt\big(|C_0'|\ge n/(s+1)^2\big) \le \exp\big({-}an/(s+1)^2\big),
\]
completing the proof of Theorem~\ref{th_Lsdecay}.
\end{proof}

\section{Percolation in other lattices}\label{sec_other}

The arguments given in~\cite{ourKesten} were specific to the case of bond percolation in $\Z^2$,
since we were trying to give as simple a proof as we could
that $p_H=1/2$ in this case. However, parts of the proofs
are applicable in many other contexts. In particular, the method used in
\SectionsecS\ of~\cite{ourKesten} applies to any planar lattice, and can be extended to other contexts.
The heart of the method is a simple application of Theorem~\ref{th_sharp};
we present this in the setting of a general lattice as Lemma~\ref{l_oursharp} in the next subsection.

In fact, the method of~\cite{ourKesten} was developed in~\cite{Voronoi} in a rather different,
continuous, context,
namely random Voronoi percolation; in~\cite{Voronoi} it is shown
that the critical probability for random Voronoi percolation in the plane is $1/2$.
The arguments needed for the random Voronoi case are much
more complicated than those for lattices; we shall not even outline them here.

In order to apply Theorem~\ref{th_sharp} to deduce results about critical probabilities,
one needs an
appropriate equivalent of the Russo-Seymour-Welsh Theorem, stating essentially that
if (very large) squares may be crossed with significant probability, then the same applies to
rectangles with a fixed aspect ratio.
As in~\cite{ourKesten}, in many contexts simpler methods can be used to prove
an essentially equivalent result.
To illustrate this we give two examples, in Subsections \ref{ssite} and \ref{stri}. The first, site percolation in the square lattice,
shows that knowing the critical probability is not necessary.
The second, site percolation in the triangular lattice, shows
that the square geometry is not necessary.

\subsection{Sharp thresholds in lattices}\label{ss_sth}

In this subsection we consider percolation on lattices in $\RR^d$.
We say that $L$ is a $d$-dimensional {\em lattice graph}, or simply {\em lattice}, if $L$ is a 
connected, locally finite
graph on a vertex set $V=V(L)\subset \RR^d$
with any two vertices at distance at least some $\rho>0$,
such that there are $d$ automorphisms $\alpha_i$
of $L$ acting on $V$ by translation through linearly independent vectors ${\bf v}_i\in \RR^d$.
We work throughout with site percolation on the graph $L$: for bond percolation we may realize
the line graph of $L$ as a lattice $L'$ and work with site percolation on $L'$.
Note that in the 2-dimensional case, $L$ need not be a planar graph.

A basic property of any lattice graph is that its vertex set
$V$ has a partition into finitely many classes $V_j$ so that the automorphism group of
the graph $L$ acts transitively on each $V_j$.

We shall need the following slightly strengthened form  of Theorem~\ref{th_sharp}.

\begin{lemma}\label{l_sharp'}
Let $X$ be a finite ground set with $|X|=N$, and suppose that $\cA\subset \Pow(X)$ is increasing.
Suppose also that there is a group $G$ acting on $X$ so that every orbit of the action
of $G$ on $X$ has size at least $M$, and so that $\cA$ is a union of orbits of the induced action
of $G$ on $\Pow(X)$.
There is an absolute constant $c_1$ such that if
$\Pr_p^X(\cA)>\epsilon$, then $\Pr_q^X(\cA)>1-\epsilon$ whenever
\[
 q-p \ge c_1\frac{\log(1/(2\epsilon))}{\log N}\frac{N}{M}.
\]
\end{lemma}

\begin{proof}
The proof is the same as that of Theorem~\ref{th_sharp}, i.e., of Theorem 2.1
of Friedgut and Kalai~\cite{FK}. Following the
proof in \cite{FK} step by step, the only modification
is that having found one variable with influence
at least $x$, one concludes that the sum of the influences of all variables is at least $Mx$,
rather than at least $Nx$.
\end{proof}

For notational convenience, we state the following result only in the $2$-dimensional case.
In $d$-dimensions corresponding results concerning paths from one face of a hypercuboid to the
opposite face, or surfaces separating one face from the opposite face, can be proved in exactly
the same way.

We work with the probability measure
$\Pr_p=\Pr_p^{L,\mathrm{site}}$ in which each vertex of $L$ is open independently with probability
$p$. An {\em open path} is a path in $L$ all of whose vertices are open.
If $L$ is a $2$-dimensional lattice and $R\subset \RR^2$ is a rectangle, then we 
write $H(R)=H_L(R)$ for the event that $R$ has a {\em horizontal open crossing}, i.e.,
that there is a path in $L$ consisting of open vertices of $R$
joining vertices $v_1$ and $v_2$, where $v_1$ is incident with an edge of $L$
that meets the left-hand side of $R$, and $v_2$ with an edge that meets the right-hand side of $R$.
In fact, for the application below the precise definition of $H(R)$ (i.e., how we deal with vertices
near the boundary of $R$) will not matter -- the statement of our lemma will not be affected
if the dimensions of the rectangles involved are altered by $O(1)$.

In this section, all our rectangles have a fixed orientation, which we take without loss
of generality to be parallel to the coordinate axes. We also suppose that the origin is a lattice point.
Note that $\Pr_p(H(R))$ may depend not just on the dimensions of $R$, but
also on its position with respect to $L$; we do not assume that the corners of our rectangles are lattice points.
In the case $L=\Z^2$, this assumption might be natural, but it would make no difference -- the statement
of the lemma is unaffected if we round the coordinates to integers.

\begin{lemma}\label{l_oursharp}
Let $L$ be a $2$-dimensional lattice graph.
Let $0<p_1<p_2<1$, $\epsilon>0$, and positive real numbers $x_1>x_2$, $y_1<y_2$ be fixed.
There is an $n_0$ such that if $n\ge n_0$ and
$R$ is an $x_1 n$ by $y_1 n$ rectangle
for which $\Pr_{p_1}(H_L(R))\ge \epsilon$,
then
$\Pr_{p_2}(H_L(R'))\ge 1-\epsilon$
for any $x_2 n$ by $y_2 n$ rectangle $R'$.
\end{lemma}

\begin{proof}
The argument is essentially the same as in~\cite{ourKesten}; we write it out for completeness.
Throughout this proof $n_0$ will be a large constant to be chosen later, depending on all
the parameters in the statement of the lemma.
 
Let ${\bf v}_1$ and ${\bf v}_2$ be two linearly independent vectors such that translations
of $\RR^2$ through ${\bf v}_i$ induce automorphisms of $L$, and let $F$ be the corresponding fundamental
region of $L$, i.e., the parallelogram with vertices $0$, ${\bf v}_1$, ${\bf v}_2$ and
${\bf v}_1+{\bf v}_2$. Note that $F$ has diameter $D=O(1)$, where the constant depends only on $L$,
and $F$ contains $\Theta(1)$ vertices of $L$.

Suppressing the dependence on $L$, suppose that
$\Pr_{p_1}(H(R))\ge \epsilon$ for an $x_1 n$ by $y_1 n$ rectangle $R$ with $n\ge n_0$.
We may find points ${\bf w}_1$ and ${\bf w}_2$, each of the form $a_1{\bf v}_1+a_2{\bf v}_2$, $a_i\in \Z$,
within distance $D$ of $((x_1+1)n,0)$
and $(0,(y_2+1)n)$, respectively. Let $F'$ be the parallelogram with vertices $0$, 
${\bf w}_1$, ${\bf w}_2$ and ${\bf w}_1+{\bf w}_2$. Then we may assume that $R$ lies within $F'$,
and indeed that $R$ does not come closer than a distance $n/3$ to the boundary of $F'$. To see this,
note that $\Pr_p(H(R))$ is unchanged if we translate $R$ through a vector ${\bf v}_i$, $i=1,2$.

Let $T$ be the graph obtained from $L$ by quotienting by (the automorphisms whose action corresponds to)
translations of $\RR^2$ through ${\bf w}_1$ and ${\bf w}_2$.
Then $T$ is a graph with $\Theta(n^2)$ vertices, where the implicit constants depend on $L$, $x_1$ and $y_2$,
and $T$ is `locally isomorphic' to $L$.
In particular, for rectangles $R'$ too small to `wrap around' $T$,
which are the only rectangles we shall consider, each rectangle $R'$ in
$L$ corresponds to a rectangle in $T$, and the induced subgraphs of $L$ and $T$ are isomorphic.

We write $\Pr_p^T$ for the probability measure in which each vertex of $T$ is open with probability
$p$, independently of all other vertices. From the remark above, 
there is an $x_1 n$ by $y_1 n$ rectangle $R$ in $T$ such that
\[
 \Pr_{p_1}^T(H(R)) = \Pr_{p_1}^{L}(H(R)) \ge \epsilon.
\]
Let $E$ be the event that there is {\em some} $x_1 n$ by $y_1 n$ rectangle $R'$ in $T$ for which
$H(R')$ holds. Then
\[
 \Pr_{p_1}^T(E) \ge \Pr_{p_1}^T(H(R)) \ge \epsilon.
\]
The event $E$ is increasing and symmetric in the sense of Lemma~\ref{l_sharp'};
translations of $T$ through the vectors ${\bf v}_i$ preserve $E$, and such translations
map any vertex of $T$ to a vertex in one given fundamental region. Thus the action of
the group generated by these translations on $T$ has $O(1)$ orbits, each of size at least $S=cn^2$,
where $c$ depends on $L$, $x_1$ and $y_2$.
We claim that for any constant $\eta<1$ we have
\[
 \Pr_{p_2}^T(E) \ge 1-\eta,
\]
provided that $n_0$ is chosen large enough, which we shall assume from now on.
Indeed, writing $N=|T|=\Theta(n^2)$, then as $N/S$ is bounded, by
Lemma~\ref{l_sharp'}
it suffices to choose $n_0$ large enough that for $n\ge n_0$ we have
$\log N=2\log n+O(1)$ larger than a certain constant depending on $\eta$ and
the parameters of the lemma.

Let $R'$ be any $x_2 n$ by $y_2 n$ rectangle in $T$. Note that
$x_1>x_2$ and $y_1<y_2$, so $R'$ is `shorter and fatter' than $R$.
It follows that if $n$ is large enough, the torus $T$
can be covered by a bounded number $M$ of translates $R_i$ of $R'$ through vectors of the form
$a_1{\bf v}_1+a_2{\bf v}_2$, $a_i\in \Z$, in such a way that any $x_1 n$ by $y_1 n$
rectangle $R$ in $T$ crosses some $R_i$ horizontally, meaning that the intersection of $R$
and $R_i$ is an $x_2 n$ by $y_1 n$ rectangle. It follows that
any horizontal open crossing of $R$ contains a horizontal open crossing of $R_i$.
Hence, if $E$ holds, then so does one of the events $E_i=H(R_i)$, so $E^c\supset \cap_i E_i^c$.

The events $E_i$ are increasing. Hence, by Lemma~\ref{l_Harris},
for each $i$ the decreasing event $E_i^c$ is positively correlated with the
decreasing event $\bigcap_{j<i} E_j^c$, and
\[
 \Pr_{p_2}^T(E^c) \ge \Pr_{p_2}^T\left(\bigcap_{i=1}^M E_i^c\right)
 \ge \prod_{i=1}^M\Pr_{p_2}^T(E_i^c) = \Pr_{p_2}^T(H(R')^c)^M.
\]
For the last step we use the fact that
the subgraph of $T$ induced by each $R_i$ is isomorphic to that induced by $R'$.
Thus,
\[
 \Pr_{p_2}^T(H(R')^c) \le \Pr_{p_2}^T(E^c)^{1/M} \le \eta^{1/M} = \epsilon,
\]
if we choose $\eta$ appropriately.
Using the local isomorphism between $L$ and $T$, we have
\[
 \Pr_{p_2}^L(H(R')) = \Pr_{p_2}^T(H(R')) \ge 1-\epsilon,
\]
as required.
\end{proof}

\subsection{Site percolation in the square lattice}\label{ssite}

For this subsection, let $\Ls=\Z^2$ be the planar square lattice viewed as a graph as in Section~\ref{bdecay}, and
let $\Lsd$ be the (non-planar) graph with vertex set $\Z^2$ in which two vertices are adjacent
if they are at Euclidean distance $1$ or $\sqrt{2}$.
We consider the probability measure $\Pr_p$ in which each vertex $v\in \Z^2$
is open with probability $p$, independently of the other vertices.
Note that we are considering two notions of site percolation involving the same probability measure.
For $L=\Ls$ or $\Lsd$, the open cluster $C_0(L)$ containing the origin is the set of open vertices
that may be reached from the origin by a path in the graph $L$ all of whose vertices are open.
As before, for a rectangle $R$ with integer coordinates, we write $H_L(R)$ for the
event that $R$ has a horizontal open crossing in $L$, and $V_L(R)$ for the event
that $R$ has a vertical open crossing.

The lattices $\Ls$ and $\Lsd$ are {\em dual} in a sense illustrated by the
following lemma.

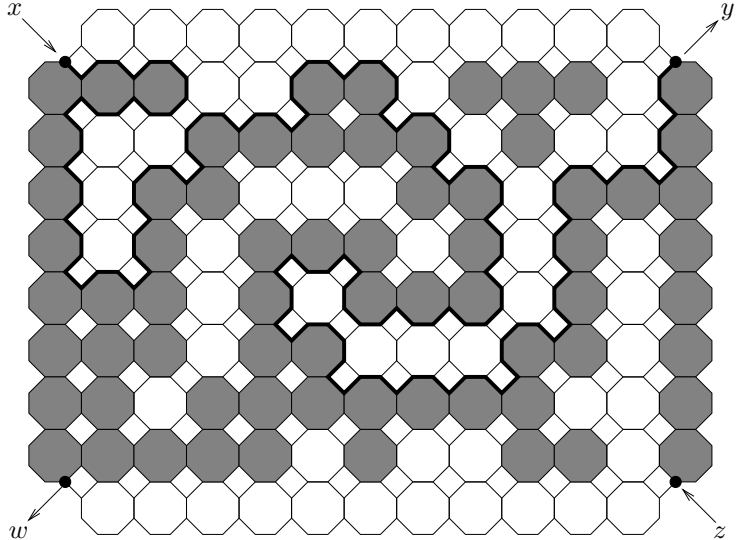
\begin{figure}[htb]
 \centering
 \input{sitedual.pstex_t}
 %\[\epsfig{file=sitedual.eps,height=2.5in}\]
 \caption{A rectangle $R$ in $\Z^2$ with each vertex drawn as an octagon, with an additional
row/column of vertices on each side. `Black' (shaded) octagons are open.
Either there is a black
path from left to right, or a white path (which may use the squares) from top to bottom.
The path $W$ entering at $x$ is shown by thick lines. As
$W$ leaves at $y$, $H_{\Ls}(R)$ holds.}
\label{fig_sitedual}
\end{figure}

\begin{lemma}\label{l_hlvl}
Let $L$ be one of $\Ls$ and $\Lsd$, let $L^\star$ be the other, and let $R$ be a
rectangle with integer coordinates. Whatever the states of the vertices in $R$,
either there is an open $L$-path crossing
$R$ from left to right, or a closed $L^\star$-path crossing $R$ from top to bottom, but
not both.
In particular,
\begin{equation}\label{hlvl}
 \Pr_p(H_L(R)) + \Pr_{1-p}(V_{L^\star}(R)) =1.
\end{equation}
\end{lemma}
\begin{proof}
Without loss of generality, we may take $L=\Ls$.
Consider the partial tiling of the plane by octagons and squares shown in Figure~\ref{fig_sitedual}:
we take one octagon for each vertex $v$ of $R$, coloured black if $v$ is open and white
if $v$ is closed, plus additional black
octagons to the left and right of $R$ and white octagons above and below $R$ as shown.
All squares are white.
Let $G$ be the graph formed by taking those edges of octagons/squares that separate a
black region from a (bounded) white one, with the endpoints
of these edges as the vertices. Then every vertex of $G$ has degree exactly $2$
except for the four vertices $x$, $y$, $z$ and $w$, which have degree $1$.
Thus the component of $G$ containing $x$ is a path $W$, ending either at $y$ or at $w$;
the path $W$ cannot end at $z$ as, walking along $W$ from $x$, one always has a black region
on the right and a white one on the left.

The black octagons on the right of $W$ correspond to an $\Ls$-connected set of sites, while
the white octagons on the left correspond to an $\Lsd$-connected set of sites.
This, if $W$ ends at $y$, as shown, there is an open $\Ls$-path from the left of $R$ to the right.
If $W$ ends at $w$, there is a closed $\Lsd$-path from the top of $R$ to the bottom.
We cannot have both crossings, as otherwise $K_5$ could be drawn in the plane.
\end{proof}

The values of the critical probabilities $p_H(L,\mathrm{site})$, $L=\Ls,\Lsd$, are not known.
A special case of the general result of Menshikov~\cite{Menshikov} (see also~\cite{MMS,Grimmett})
implies that for $L=\Ls$ or $L=\Lsd$ there is
exponential decay of the radius of $C_0(L)$ below $p_H(L,\mathrm{site})$, and hence that
$p_T(L,\mathrm{site})=p_H(L,\mathrm{site})$. As noted in the introduction, it follows
from the results of Kesten~\cite{Kesten81} or
Aizenman and Newman~\cite{AN} (see also~\cite{Grimmett})
that there is exponential decay of $|C_0(L)|$.
We give a new proof of the latter, stronger result.

\begin{theorem}\label{th_sitesq}
Let $L=\Ls$ or $\Lsd$. For any $p<p_H(L,\mathrm{site})$, there is a constant $a=a(p,L)>0$ such that
$\Pr_p(|C_0(L)|\ge n)\le \exp(-an)$ holds for all $n\ge 0$.
\end{theorem}

In proving Theorem \ref{th_sitesq} we shall make use of the
following more general version of \LemmaX{} of~\cite{ourKesten}.
When there is no danger of ambiguity, we write $H(R)$ for $H_L(R)$ and $V(R)$
for $V_L(R)$.

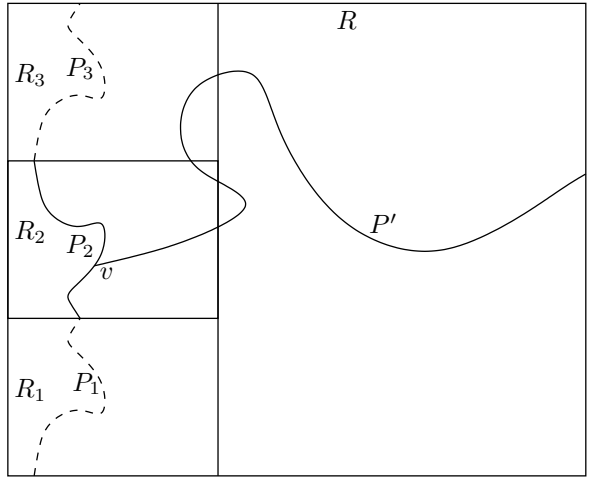
\begin{figure}[htb]
\centering
\input{l_X3.pstex_t}
\caption{The rectangles $R_i$ and rectangle $R$ for $k=3$: the solid paths indicate that $X_2$ holds.}
\label{fig_l_Xk}
\end{figure}

\begin{lemma}\label{l_Xk}
Let $L=\Ls$ or $\Lsd$, and let $k$, $r$, $s$ and $t>r$ be positive integers.
Set $R_i=[0,r]\times [(i-1)s,is]$ for $i=1,2,\ldots,k$, and let $R=[0,t]\times [0,ks]$.
Let $X_i$ be the event that there is an open vertical crossing
of $R_i$ joined by an open path in $R$ to the right-hand side of $R$.
Then for some $i$ we have
\[
 \Pr_p(X_i) \ge \Pr_p(H(R))\Pr_p(V(R_1))/k.
\]
\end{lemma}

\begin{proof}
The proof is almost exactly the same as that of \LemmaX{} of~\cite{ourKesten}.
If $V(R_i)$ holds, we can define a left-most vertical crossing $\LP(R_i)$ of $R_i$
in such a way that the event $\LP(R_i)=P_i$ does not depend on the states of vertices of $R_i$ to the right
of $P_i$. (This is illustrated rotated in Figure \ref{fig_sitedual}: there is a horizontal open
crossing $P$ consisting of sites next to the path $W$. Finding $W$ step by step, we only ever
examine vertices adjacent to $W$, so no vertex below $P$ has been examined.)

For a fixed $i$, if $V(R_i)$ holds and $\LP(R_i)=P_i$, define $P$ to be the vertical (but not necessarily
open) crossing of $[0,r]\times [0,sk]$ obtained by reflecting $P_i$ in the horizontal lines $y=js$,
as shown in Figure~\ref{fig_l_Xk}. Also, let $P_j$, $1\le j\le k$, be the sub-paths of $P$ crossing
each $R_j$. Note that the event that $P$ takes a particular value
is independent of the states of the vertices to the right of $P$.

With (unconditional) probability $\Pr_p(H(R))$ there is a horizontal open crossing $P_H$ of $R$. Any such
crossing must cross $P$; indeed, $P$ and $P_H$ share a vertex unless $L=\Lsd$ and the paths
cross diagonally within a grid square. It follows that $P_H$
contains a sub-path $P'$ with the following properties:
every vertex of $P'$ lies strictly
to the right of $P$ and is open, $P'$ starts at a vertex adjacent to a vertex $v$ of $P$,
and $P'$ ends at a vertex on the right hand side of $R$; see Figure~\ref{fig_l_Xk}.
Let $Y_j(P)$ be the event that such a $P'$ exists
with $v$ lying on $P_j$. Then we have
\begin{equation}\label{cs1}
 \sum_{j=1}^k \Pr_p(Y_j(P)) \ge \Pr_p(H(R)).
\end{equation}
Now $Y_j(P)$ depends only on the states of the vertices to the right of $P$.
For any possible value $P_i$ of $\LP(R_i)$, defining $P$ and $P_j$ as above,
the event $\LP(R_j)=P_j$ is independent of the states of vertices to the right of the path $P$.
Thus,
\[
 \Pr_p(Y_j(P)\mid \LP(R_j)=P_j) = \Pr_p(Y_j(P)),
\]
and, from \eqref{cs1},
\[
 \sum_{j=1}^k \Pr_p(Y_j(P) \mid \LP(R_j)=P_j) \ge \Pr_p(H(R)).
\]
If $Y_j(P)$ holds and $\LP(R_j)=P_j$, then $X_j$ holds (see Figure~\ref{fig_l_Xk}). Thus,
\[
 \sum_{j=1}^k \Pr_p(X_j \mid \LP(R_j)=P_j) \ge \Pr_p(H(R)).
\]
In other words,
\[
 \sum_{j=1}^k \frac{\Pr_p(X_j\hbox{ holds and }\LP(R_j)=P_j)}
 {\Pr_p(\LP(R_j)=P_j)} \ge \Pr_p(H(R)).
\]
Recalling the definition of the paths $P_j$, we have
\[
 \Pr_p(\LP(R_j)=P_j)=\Pr_p(\LP(R_1)=P_1),
\]
so
\begin{equation}\label{eq}
 \sum_{j=1}^k \Pr_p(X_j\hbox{ holds and }\LP(R_j)=P_j)
 \ge \Pr_p(H(R)) \Pr_p(\LP(R_1)=P_1).
\end{equation}

So far, $P_1$ was fixed.
As $P_1$ runs over all possible values of $\LP(R_1)$, each $P_j$ runs over all possible values of $\LP(R_j)$.
Summing \eqref{eq} over $P_1$, as $V(R_j)$ is the disjoint union of the events
that $\LP(R_j)$ takes each possible value, it follows that
\[
 \sum_{j=1}^k \Pr_p(X_j) \ge \Pr_p(H(R)) \Pr_p(V(R_1)),
\]
and the result follows.
\end{proof}

As in~\cite{ourKesten},
we obtain an immediate corollary concerning long thin rectangles, provided we know that
certain crossings of squares exist with significant probability.
We write $R_{m,n}$ for the $m$ by $n$ rectangle
$[0,m]\times [0,n]$, and $H(R_{m,n})$ for the event that this rectangle
has a horizontal open crossing in the lattice under consideration.

\begin{corollary}\label{c_longss}
Let $c>0$ and integers $\rho$, $k\ge 2$ be given. There is a constant $c'=c'(c,k,\rho)>0$ such
that if $L=\Ls$ or $\Lsd$, and $\Pr_p(H(R_{s,s})),\Pr_p(H(R_{ks,ks}))\ge c$,
then $\Pr_p(H(R_{\rho ks,ks}))\ge c'$.
\end{corollary}

\begin{proof}
Let $h_{m,n}=\Pr_p(H(R_{m,n}))$, so $h_{s,s}$, $h_{ks,ks}\ge c$ by assumption.
We claim that for $m>s$ we have
\begin{equation}\label{hms}
 h_{2m-s,ks}\ge h_{m,ks}^2 c^3/k^2.
\end{equation}
Applying \eqref{hms} this repeatedly, the result follows.

As in~\cite{ourKesten}, the inequality \eqref{hms}
is an immediate consequence of Lemma~\ref{l_Xk} and Harris' Lemma.
To see this, choose an $i$, $1\le i\le k$, for which Lemma~\ref{l_Xk} holds with $r=s$, $t=m$,
and consider the rectangles $R=[0,m]\times [0,ks]$, $R'=[s-m,s]\times [0,ks]$ and the square
$S=[0,s]\times [(i-1)s,is]$ in their intersection. Note that the square $S$ plays the role of the
rectangle $R_i$ in Lemma~\ref{l_Xk} for the parameters ($r=s$, $t=m$) we have used.

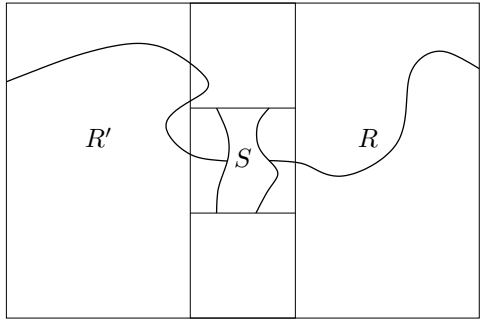
\begin{figure}[htb]
 \centering
 \input{RSRk.pstex_t}
\caption{The overlapping rectangles $R$ and $R'$ with the square $S$ in their intersection. The paths drawn show
that $X_i$ holds for $R$, as well as the reflected equivalent $E_2$ for $R'$.
If $H(S)$ also holds, then so does $H(R\cup R')$.}
\label{fig_RSRk}
\end{figure}

Let us write $E_1$ for the event $X_i$ defined in Lemma~\ref{l_Xk}, which depends on the vertices in $R$,
and let $E_2$ be the corresponding event for $R'$, defined by reflecting in the line $x=s/2$;
see Figure~\ref{fig_RSRk}.
Finally, let $E_3$ be the event $H(S)$. Note that if $E_1$, $E_2$ and $E_3$ all hold, then
$H(R\cup R')$ holds, using only the fact that horizontal and vertical crossings of $S$ must cross.
By Lemma~\ref{l_Xk} and our choice of $i$ we have 
\[
 \Pr_p(E_1)\ge \Pr_p(H(R))\Pr_p(V(S))/k = h_{m,ks} \Pr_p(V(S))/k.
\]
By symmetry, $\Pr_p(E_2)=\Pr_p(E_1)$. As $E_1$, $E_2$ and $E_3$ are increasing events, by Lemma~\ref{l_Harris} we have
\begin{eqnarray*}
 \Pr_p(H(R\cup R')) &\ge& \Pr_p(E_1)\Pr_p(E_2)\Pr_p(E_3) \\
 &\ge& h_{m,ks}^2 \Pr_p(V(S))^2\Pr_p(H(S))/k^2.
\end{eqnarray*}
By assumption, $\Pr_p(V(S))=\Pr_p(H(S))=h_{s,s}\ge c$, so
\[
 h_{2m-s,ks} = \Pr_p(H(R\cup R')) \ge h_{m,ks}^2 c^3/k^2,
\]
completing the proof of \eqref{hms} and thus of the corollary.
\end{proof}

Using the method of \SectionsecS\ of~\cite{ourKesten}, it is easy to deduce Theorem~\ref{th_sitesq}.
The key step is to apply Lemma~\ref{l_oursharp}.

\begin{proof}[Proof of Theorem~\ref{th_sitesq}]
Let $L=\Ls$ or $\Lsd$.
It suffices to show that for any constant $p_1<p_2$, either percolation occurs in $L$
at $p_2$ (i.e., $\theta_L(p_2)>0$), or there is exponential decay of $|C_0(L)|$ at $p_1$.
Fix $p_1<p_2$, and set $p=(p_1+p_2)/2$.
Let $n_0$ be a large constant to be chosen later, depending only on $p_1$ and $p_2$.

For $i=1,2,4$, let $S_i$ be a square of side length $in_0$. From Lemma \ref{l_hlvl},
we have $\Pr_p(H_L(S_i))+ \Pr_{1-p}(H_{L^\star}(S_i)) =1$, where $\{L,L^\star\}=\{\Ls,\Lsd\}$.
It follows that either (a) there are two values of $i\in \{1,2,4\}$ for which
$\Pr_p(H_L(S_i))\ge 1/2$, or (b) there are two values for which $\Pr_{1-p}(H_{L^\star}(S_i))\ge 1/2$.

It follows from Corollary~\ref{c_longss} (applied with $c=1/2$, $\rho=10$, and $k=2$ or $k=4$)
that there is a $10n$ by $n$ rectangle $R_n$, $n\ge n_0$, such that
\begin{equation}\label{pq}
 \Pr_p(H_L(R_n))\ge c' \hbox{ or }\Pr_{1-p}(H_{L^\star}(R_n))\ge c',
\end{equation}
where $c'$ is an absolute constant not depending on our choice of $n_0$.

As $p_1<p<p_2$, for any constant $c_3<1$ it follows by Lemma~\ref{l_oursharp}
that if $n_0$ was chosen large enough, and $R'$ is a $6n$ by $2n$ rectangle
with $n$ as above, then either
\begin{equation}\label{e_p}
 \Pr_{p_2}(H_L(R'))\ge c_3
\end{equation}
or
\begin{equation}\label{e_dec}
 \Pr_{1-p_1}(H_{L^\star}(R')) \ge c_3.
\end{equation}

If~\eqref{e_dec} holds and $c_3$ is chosen large enough then,
as an open path in $L$ cannot start inside and end outside a closed cycle
in $L^\star$, we can use Lemma~\ref{l_kdepneg} exactly as in Section~\ref{bdecay}
to obtain exponential decay of the size of $C_0(L)$ in $\Pr_{p_1}$.

If~\eqref{e_p} holds and $c_3$ is chosen large enough, then $\theta_L(p_2)>0$ follows.
There are several standard arguments;
we outline a slightly less standard one, given in~\cite{ourKesten}.
Choose $c_3=p_0^{1/3}$, where $p_0$ is some constant for
which Lemma~\ref{l_1dep} holds. For a $6n$ by $2n$ rectangle $R$, let $G(R)$ be the event
that $H(R)$, $V(S_1)$ and $V(S_2)$ all hold, where the $S_i$ 
are the two $2n$ by $2n$ `end' squares
of $R$. Note that $\Pr_{p_2}(V(S_i))=\Pr_{p_2}(H(S_i))\ge \Pr_{p_2}(H(R))\ge c_3$.
Thus, by Lemma~\ref{l_Harris}, $\Pr_{p_2}(G(R))\ge c_3^{\,3}=p_0$.
We define a $1$-dependent bond percolation measure $\Prt$
on $\Z^2$ by declaring the edge from $(a,b)$ to $(a+1,b)$ to be open in $\Prt$
if $G(R)$ holds in $\Pr_{p_2}$
for the $6n$ by $2n$ rectangle with bottom left corner $(2an,2bn)$. The definition for vertical edges
is analogous. By Lemma~\ref{l_1dep} we have percolation in $\Prt$. The definition of the event $G(R)$
ensures that for any open path $P$ in $\Prt$ there is a corresponding open path $P'$ in $L$.
When $P$ is infinite, so is $P'$, so site percolation occurs in $L$ in the probability
measure $\Pr_{p_2}$, i.e., $\theta_L(p_2)>0$.
\end{proof}

Theorem \ref{th_sitesq} certainly implies that $p_T(L)=p_H(L)$ for $L=\Ls$ or $\Lsd$.
Together with an intermediate step \eqref{pq} in the proof above, it also implies
the well-known result relating $p_H(\Ls)$ and $p_H(\Lsd)$.

\begin{corollary}
For site percolation we have $p_H(\Ls)+p_H(\Lsd)=1$.
\end{corollary}
\begin{proof}
Suppose first that $p_H(\Ls)+p_H(\Lsd)>1$.
Then there is a $p$ with $p<p_H(\Ls)$ and $1-p<p_H(\Lsd)$. By
Theorem \ref{th_sitesq}, we have
exponential decay of $|C_0(\Ls)|$ in $\Pr_p$ and of $|C_0(\Lsd)|$ in $\Pr_{1-p}$.
Thus the $\Pr_p$-probability that a large square has either a horizontal open $\Ls$-crossing
or a vertical closed $\Lsd$-crossing tends to zero, contradicting Lemma \ref{l_hlvl}.

It remains to show that $p_H(\Ls)+p_H(\Lsd)\ge 1$,
which is analogous to Harris' Theorem for bond percolation.
To show this, we shall prove that any $p$ we have
\begin{equation}\label{one0}
 \theta_\Ls(p)=0 \hbox{ or } \theta_\Lsdsub(1-p)=0.
\end{equation}
This follows from \eqref{pq} in a standard
way, analogous to the proof of Harris' Theorem, \TheoremHarris, in \cite{ourKesten}.
Indeed, from \eqref{pq} there is a sequence $n_i$ with $n_{i+1}\ge 4n_i$ such that
for each $i$, either $\Pr_p(H_{\Ls}(R_{n_i}))\ge c'$, or $\Pr_{1-p}(H_\Lsdsub(R_{n_i}))\ge c'$.
Passing to a subsequence $m_i$, we may assume that one case always holds. If the first case
holds, then we may construct annuli $A_i$ as in Figure \ref{fig_ann}
with inner and outer radii $m_i$ and $3m_i$, so that the $A_i$ are disjoint, and each surrounds 
the origin. By Lemma \ref{l_Harris}, each $A_i$ contains an open $\Ls$-cycle surrounding
the origin with probability at least $(c')^4$. Hence, with probability $1$ some $A_i$ 
contains such a cycle, and it follows that $\theta_\Lsdsub(1-p)=0$. Similarly, in the other
case $\theta_{\Ls}(p)=0$, proving \eqref{one0}. As noted above, $p_H(\Ls)+p_H(\Lsd)=1$ follows.

Also, we have shown that $\theta_L(p_H(L))=0$ for at least one of $\Ls$ and $\Lsd$.
\end{proof}

Let us remark that, as pointed out by Professor Ronald Meester and described in \cite{BRnote},
one can use a sharp-threshold of Russo~\cite{Russo01} in place of Lemma~\ref{l_oursharp}.

\subsection{Site percolation in the triangular lattice}\label{stri}

In this subsection we consider the equilateral triangular lattice $\Lt$ with edge length $1$.
We shall take the origin and the point $(0,1)$ on the $y$-axis to be vertices of $\Lt$.
Each vertex of $\Lt$ will be open independently
with probability $p$; we write $\Pr_p$ for this site percolation measure. As usual, $\Lt$
will be viewed as a graph, in which vertices at distance $1$ are adjacent.

It is well-known that $p_H(\Lt)=1/2$. Indeed, the following result
is another special case of the general results mentioned in the introduction.

\begin{theorem}\label{th_tri}
In the triangular lattice $\Lt$, if $p>1/2$ then $\theta(p)>0$. If $p<1/2$, then there
is a constant $a=a(p)>0$ such that $\Pr_p\big(|C_0(\Lt)|\ge n\big)\le \exp(-an)$ holds for all $n\ge 0$.
\end{theorem}

The arguments will be very similar to those in the previous sections, so we only sketch the details.

Although the natural equivalent of \Lemmahalf{} in~\cite{ourKesten} (i.e., the standard starting point
that the crossing probability
for a square is $1/2$ in $p=1/2$ bond percolation on $\Z^2$) applies to a parallelogram with a $60$ degree angle,
we shall work with rectangles; parallelograms do not seem to fit together in the
way required for the equivalent of Lemma~\ref{l_Xk}. Also, while a symmetry argument shows that
the crossing probability for a suitably oriented
square is $1/2$ at $p=1/2$, this works only for certain
orientations. These orientations will not be consistent with the symmetry required in Lemma~\ref{l_Xk}.

Unlike in previous
sections, the rectangles we consider will often not be aligned with the coordinate axes.
Given a non-square rectangle, we define {\em long} and {\em short} crossings of $R$ in the obvious way,
and write $L(R)$ and $S(R)$ respectively for the events that $R$ has a long open crossing or
a short open crossing.

As the neighbourhood of a vertex of $\Lt$ is connected,
if $C$ is a finite open cluster in $\Lt$,
then its vertex boundary contains a closed cycle $S$ surrounding $C$. Also, if a path in $\Lt$ starts inside and ends
outside a cycle, then the path and cycle share a vertex. It follows that if $R$ is not too small
(say both sides have length at least two), then $R$ has a long open crossing if and only if $R$
does not have a short {\em closed} crossing. Hence,
\[
 \Pr_p(L(R))+\Pr_{1-p}(S(R))=1.
\]
In particular,
\begin{equation}\label{selfdual}
 \Pr_{1/2}(L(R))+\Pr_{1/2}(S(R))=1.
\end{equation}

Most of the work needed to prove Theorem~\ref{th_tri} is contained in the following lemma.
Working in $\Z^2$, we took our rectangles to be aligned with the coordinate axes.
Here, we do not specify the orientation of the rectangle $R$.

\begin{lemma}\label{l_tri}
There is an absolute constant $c>0$ such that
for any $n_0$ there is an $n\ge n_0$ and a $6n$ by $n$ rectangle $R$ with
\begin{equation}\label{tric}
 \Pr_{1/2}(L(R))\ge c.
\end{equation}
\end{lemma}

\begin{proof}
The idea is to use an equivalent of Lemma~\ref{l_Xk} for $\Lt$. 
In fact, we have written the proof of Lemma~\ref{l_Xk} so that it goes through unchanged
for $L=\Lt$, noting that the lines $y=is$
that we reflect in are symmetry axes of $\Lt$.

In order to use an argument similar to that of Corollary~\ref{c_longss} to deduce Lemma~\ref{l_tri},
we need
as a starting point that certain crossing probabilities of rectangles are not too small.

Consider a fixed integer $s$, and rectangles of the form $[a,b]\times [0,s]$, where $a,b$, $b-a>2$
are integer multiples of $\sqrt{3}/2$.
If $R$ and $R'$ are two rectangles of this form with $R\subset R'$, and $R'$ is obtained
by extending $R$ horizontally by a distance of $\sqrt{3}/2$, then $R'$ contains
one extra column of lattice points.
\begin{figure}[htb]
 \[\epsfig{file=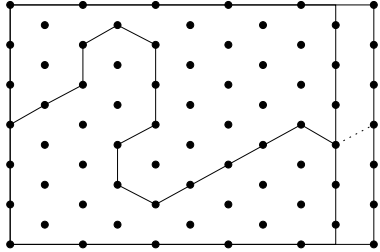,height=1.3in}\]
 \caption{A rectangle $R'$ extending a rectangle $R$ by one column of lattice points, a path crossing
$R$ horizontally (solid lines), and an extension to a crossing of $R'$ (dashed line).}
\label{fig_RRp}
\end{figure}
 As $R'$ extends $R$ horizontally, we have
$\Pr_{1/2}(H(R'))\le \Pr_{1/2}(H(R))$. However, we also have
\begin{equation}\label{err}
 \Pr_{1/2}(H(R'))\ge \Pr_{1/2}(H(R))/2.
\end{equation}
Indeed, $H(R)$ depends only on the states of points inside $R$, and if $R$ has an open crossing then there
is at least one point in $R'\setminus R$ which, if open, extends this crossing to an open crossing of $R$;
see Figure~\ref{fig_RRp}.

Suppose that Lemma~\ref{l_tri} does not hold
and, in particular, that it does not hold with $c=0.01$, say.
Then there is an $n_0$ such that for any $n\ge n_0$ and any
$6n$ by $n$ rectangle $R$ with any orientation, we have
\begin{equation}\label{efail}
 \Pr_{1/2}(L(R))<0.01.
\end{equation}
We claim that, for any integer $s\ge 6n_0$,
there is a real number $t(s)$ which is an integer multiple of $\sqrt{3}$, such that
\begin{equation}\label{pts}
 1/8 \le \Pr_{1/2}\big(H([0,t]\times [0,s])\big) \le 1/2
\end{equation}
holds for $t=t(s)$.
Indeed, as $t$ increases, the probability above decreases, and from the observation \eqref{err} above
it cannot decrease by more than a
factor of $4$ as $t$ increases by $\sqrt{3}$. Also, by \eqref{efail},
the probability above is at most $0.01$ for $t=6s$ and, using \eqref{selfdual},
at least $0.99$ for $t=s/6$.
Hence $s/6\le t(s)\le 6s$.

This gives us a starting point for the induction used in the proof of Corollary~\ref{c_longss}:
using~\eqref{pts} and~\eqref{selfdual}, we see that for $s=6n_0$, $R_1=[0,t(s)]\times [0,s]$
has $\Pr_{1/2}(H(R_1))$, $\Pr_{1/2}(V(R_1))\ge 1/8$. The same follows for $R_i=[0,t(s)]\times [(i-1)s,is]$,
as each $R_i$ is positioned in the same way with respect to the lattice as $R_1$.
The second ingredient of the starting point is the large rectangle $R$, for which we may take
$[0,t(40s)]\times [0,40s]$, using $k=40$ when we apply Lemma~\ref{l_Xk}.
Note that we have $t(40s)\ge 40s/6 = (40/36)6s \ge (40/36) t(s)$.
Now the proof of Corollary~\ref{c_longss} goes through as before, noting that all the rectangles
we consider have vertices that are lattice points, and that the line $x=t(s)/2$ is a symmetry
axis of $\Lt$.
\end{proof}

\begin{proof}[Proof of Theorem~\ref{th_tri}]
The method is similar to that we used
for the square lattice, so we give only an outline, emphasizing the differences.

\begin{figure}[htb]
 \[\epsfig{file=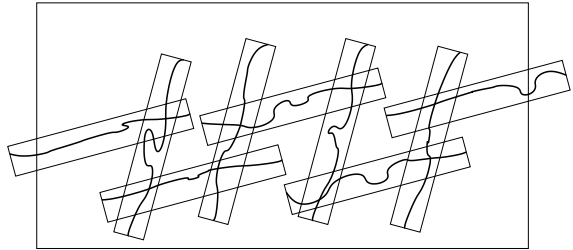,height=1.3in}\]
\caption{A path of congruent $6$ by $1$ rectangles crossing a large rectangle $R'$.}
\label{fig_rectpath}
\end{figure}

Let $n_0$ be a large constant, to be chosen below.
Let $R$ be a $6n$ by $n$ rectangle
with $n\ge n_0$ for which \eqref{tric} holds; the existence of such an $R$ is guaranteed by Lemma~\ref{l_tri}.
We first note that there is an absolute constant $c'>0$ (not depending on $n_0)$
such that if $R'$ is a $34n$ by $10n$ rectangle
with any orientation, and any position with respect to the lattice,
then $\Pr_{1/2}(L(R'))\ge c'$. To see this,
construct a path of rectangles $R_i$ inside $R'$, with each $R_i$
congruent to $R$ and placed similarly with respect to the grid,
so that long crossings of $R_i$ and $R_{i+1}$ cross, and long crossings of the first and last $R_i$
cross the opposite short sides of $R'$, as in Figure~\ref{fig_rectpath}. 
Then apply Lemma~\ref{l_Harris}, noting that the number of rectangles in the path is bounded by
some absolute constant.
(In fact, this construction is possible starting from a rectangle $R$ with any fixed aspect ratio
larger than $(1+\sqrt{3})/2$, but with a larger aspect ratio the picture is clearer.)

The rest of the proof is as for the square lattice.
Fix any $p>1/2$, and any $c_1<1$.
From the argument above and Lemma~\ref{l_oursharp}, if $n_0$ is chosen large enough, then for $n\ge n_0$
we have $\Pr_p(L(R''))\ge c_1$ for every $33n$ by $11n$ rectangle $R''$.
As in previous sections, we can apply Lemma~\ref{l_1dep} to deduce that $\theta_\Lt(p)>0$,
and Lemma~\ref{l_kdepneg} to deduce exponential decay of $|C_0(\Lt)|$ in $\Pr_{1-p}$.
\end{proof}

\section{Percolation in symmetric environments}\label{sec_non}

So far, we have considered site percolation on lattice graphs $L$. The
lattice structure was used in two ways: firstly, the notion of
percolation, or of an open crossing of a rectangle, was defined using
paths in $L$ consisting of vertices of $L$ that are open, where the
model was that the states of vertices were independent. Secondly, the
symmetry of the lattice was important, principally in the application
of Lemma~\ref{l_sharp'}. For our methods, the second use of the lattice
structure is essential, but the first is not.
Rather than write a very general version of Lemma~\ref{l_oursharp},
whose statement would be almost as long as its proof, we shall illustrate
this with two examples. In these settings
the method of Menshikov~\cite{Menshikov} does not seem to work, as the
van den Berg-Kesten inequality~\cite{vdBK} does not apply.

We start by discussing the other main ingredient of our approach,
namely, a suitable equivalent of the Russo-Seymour-Welsh (RSW) Theorem.

\subsection{A general weak RSW Theorem}

In \cite{Voronoi}, a weak version of the RSW Theorem was proved for random Voronoi percolation,
where the Voronoi cells associated to a Poisson process in the plane are coloured.
Due to the more complicated setting, the proof of this result, \TheoremVRSW\ of~\cite{Voronoi},
is rather long. However, as noted in~\cite{Voronoi}, the result holds for a wide class
of percolation models.
While weaker than the natural analogue of the RSW Theorem (whose truth is not known for random
Voronoi percolation), the result in~\cite{Voronoi} is strong enough to serve
as a key step in establishing the critical probability.

Certain properties of the crossings that arise in percolation models are rather general.
For example, in either bond percolation
in $\Z^2$ or random Voronoi percolation, horizontal and vertical open crossings of the same rectangle
must meet. Hence, such crossings of suitably arranged overlapping rectangles can be combined to form
crossings of longer rectangles. To generalize this observation, we may consider any probability
measure that assigns a {\em state}, {\em open} or {\em closed}, to each point $x$ of some set $S\subset \RR^2$.
In this setting, a horizontal open crossing of a rectangle $R\subset \RR^2$ is a (piecewise-linear)
geometric path $P\subset R$ starting at a point on the left-hand side of $R$ and ending at
a point on the right-hand side, such that
every point of $P$ is open. In the bond percolation case, we may take $S$ to be the set
of points with at least one integral coordinate; this set is exactly the subset of $\RR^2$
obtained 
when we draw the graph $\Z^2$ with straight-line segments as edges. A point of $S\setminus \Z^2$
is open if the corresponding edge is open. We may take the points of $\Z^2$
to be always open. Then crossings by open paths in the graph $\Z^2$ correspond to open paths $P$
as defined above. In the random Voronoi setting we have $S=\RR^2$, and a point of $S$ is open if it lies
in an open Voronoi cell (defined with respect to a Poisson process).

Below we shall restate \TheoremVRSW\ of~\cite{Voronoi} as Theorem~\ref{th_wRSW};
in~\cite{Voronoi}, this result was formally stated and proved only for random Voronoi percolation,
but it was noted that the proof given applies essentially without modification in a much more general
setting, which we now describe. 

Let us suppose that we have a probability measure $\Pr$ on assignments of a state,
open or closed, to each point of some subset $S$ of $\RR^2$, with the following
additional assumptions.

(i) The event
that a point, or a measurable subset, of $S$ is open is increasing in a suitable
product space, so that Lemma~\ref{l_Harris} can be applied to events such as `$R$ 
has a horizontal open crossing'.

(ii) The set-up
has the symmetries of $\Z^2$, i.e., is unchanged by translation through the vectors
$(1,0)$ and $(0,1)$, reflection in the axes, and rotation through 90 degrees
about the origin.

(iii) Disjoint regions are asymptotically independent as we `zoom out'.
To make this precise, for $R\subset \RR^2$
and $\lambda\in \RR$ let us write
$\lambda R$ for $\{\lambda x:x\in R\}$. We assume that if $R_1$ and $R_2$
are disjoint rectangles, then for any $\epsilon>0$ there is a $\lambda_0$
such that for any $\lambda>\lambda_0$ and any events $A_1$ and $A_2$ defined in terms
of the states of points in $\lambda R_1$ and $\lambda R_2$ respectively,
we have $|\Pr(A_1\cap A_2)-\Pr(A_1)\Pr(A_2)|\le \epsilon$.

(iv) Shortest paths are not too long: there is a constant $C$ such that,
for any fixed rectangle $R$, the probability that $H(\lambda R)$ holds but the shortest
open path $P$ crossing $\lambda R$ has length at least $\lambda^C$ tends to
zero as $\lambda\to\infty$.

Note that all these assumptions hold in the random Voronoi setting (see~\cite{Voronoi}). Also, they hold
for bond percolation in $\Z^2$; for example, for (iv) note that any shortest
open crossing of an $m$ by $n$ rectangle uses each vertex
at most once and hence has length at most $(m+1)(n+1)$.
We shall describe other settings in which the assumptions above hold
in the subsequent subsections. We write $R_{m,n}$ for the $m$ by $n$ rectangle $[0,m]\times [0,n]$.

\begin{theorem}\label{th_wRSW}
Let $c>0$ and $\rho>1$ be given. Under the assumptions above,
if $\Pr(H(R_{n,n}))\ge c$ for all large enough $n$, then there is a $c'>0$ such that
for any $n_0$ there is an $n>n_0$ with $\Pr(H(R_{\rho n,n}))>c'$.
\end{theorem}

\begin{proof}
\TheoremVRSW\ of~\cite{Voronoi} states
that, for random Voronoi percolation, for any $\rho>1$,
$\liminf\Pr(H(R_{n,n}))>0$ implies $\limsup\Pr(H(R_{\rho n,n}))>0$.
As noted in~\cite{Voronoi}, the proof uses only the assumptions
above, so the same conclusion holds in the setting here.
This is exactly our conclusion here.
\end{proof}

\subsection{Dependent bond percolation}

In this section we shall show that our methods can be applied
to dependent percolation as well. Our example is
a particular model of bond percolation on $\Z^2$,
where the states of the edges are not independent.
(As far as we are aware, this model has not been previously considered.)
Let $\hZs$ consist of the points $x=(a,b)$ with $2a,2b\in \Z$.
Our underlying probability space will consist of independent identically distributed
$\{-1,+1\}$-valued random variables $v_x$, $x\in \hZs$, with $\Pr(v_x=+1)=p$.
Let $w$ be a function from $\hZs$ to $\Z$ with the following properties:
$w(a,b)\ge 0$ for all $(a,b)$, $w$ has finite support, $w(0,0)$ is odd,
$w(a,b)$ is even unless $a=b=0$, and $w(a,b)=w(b,a)=w(-a,b)$
for all $(a,b)\in \hZs$, so $w$ has the rotational and reflectional symmetries of $\Z^2$.
We assign states to the edges of $\Z^2$ as follows: an edge $e$
of $\Z^2$ has a midpoint $m(e)\in \hZs$. Let $e$ be open
if
\begin{equation}\label{sdef}
 \sum_{x\in\hZs} w(x)v_{m(e)+x}>0.
\end{equation}
Note that the sum above is always odd, and that if $p=1/2$ then $e$ is open
with probability $1/2$.

Let us write $\ppw$ for the probability measure defined above.
As before, we write $C_0$ for the open cluster containing the origin, i.e.,
the set of vertices of $\Z^2$ connected to $(0,0)$ by a path of open edges.
We write $\thpw$ for $\ppw(|C_0|=\infty)$.
Our next result shows that the Harris-Kesten result for (independent) bond percolation in $\Z^2$
extends to this particular locally-dependent setting.

\begin{theorem}\label{th_Ldep}
Let $w:\hZs\to \Z$ satisfy the conditions above.
If $p>1/2$, then $\thpw>0$. If $p<1/2$, then there is a constant $a=a(w,p)>0$ such that
$\ppw(|C_0|\ge n)\le \exp(-an)$ for all $n>0$.
\end{theorem}

We outline the proof, which is very similar to the proof of the Harris-Kesten
Theorem given in \cite{ourKesten} together with the proof of Theorem~\ref{th_Lsdecay};
note that these results are exactly the special case when $w=0$ except at the origin.

\begin{proof}[Outline proof of Theorem~\ref{th_Ldep}.]
As usual, given a rectangle with integer coordinates we write $H(R)$ ($V(R)$)
for the events that $R$ has a horizontal (vertical) crossing by open edges.
Let $L^\star$ be the dual lattice to $L=\Z^2$; we may realize $L^\star$ so that the dual
edge $e^\star$ of each edge $e$ of $L$ has the same midpoint as $e$.
As in the independent case (see \LemmaHalf{} of \cite{ourKesten}),
taking the state of $e^\star$ to be the same as the state of $e$,
$R=[a,b]\times [c,d]$ has a horizontal open crossing if and only if the corresponding dual
rectangle $R'=[a+1/2,b-1/2]\times [c-1/2,d+1/2]$ has no closed vertical crossing;
indeed, the probability measure is irrelevant to this observation.
In our set-up, the state of $e^\star$ is also defined by \eqref{sdef}.
Hence, $e^\star$ is closed if and only if
\[
 \sum_{x\in\hZs} w(x)(-v_{m(e^\star)+x})>0,
\]
and the distribution of closed edges in $\ppw$ is exactly the distribution of 
open edges in $\pqw$. Taking $R$ to be an $n+1$ by $n$ rectangle and using
the isomorphism between $L=\Z^2$ and its dual that rotates $R'$ onto $R$, it follows that
$\ppw(H(R))+\pqw(H(R))=1$. In particular, $\phw(H(R))=1/2$, as in the independent case.

Writing $R_{m,n}$ for an $m$ by $n$ rectangle, we have
\begin{equation}\label{sdfs}
 \phw(H(R_{n,n}))\ge\phw(H(R_{n+1,n}))= 1/2.
\end{equation}
Our set-up satisfies all the conditions of Theorem~\ref{th_wRSW}: we define $S\subset\RR^2$
and the states of points of $S$ exactly as in the independent case discussed above.
From \eqref{sdef}, the event that a bond is open is increasing in the product probability space
defined by the $v_x$, and condition (i) follows. Condition (ii) follows from our symmetry assumptions on $w$,
and (iv) is immediate as for independent bond percolation. Finally, (iii) follows from the assumption that
$w$ has finite support -- indeed, for some constant $D$
we obtain complete independence of regions separated by a distance
of at least $D$.

Using Theorem~\ref{th_wRSW} and~\eqref{sdfs},
there is a $c'>0$ such that there are arbitrarily large $n$ with $\phw(R_{10n,n})>c'$.
Fix $p>1/2$.
We claim that for any $c''<1$, there are arbitrarily large $n$ with $\ppw(R_{6n,2n})> c''$.
This follows from Lemma~\ref{l_sharp'} in essentially the same way as
Lemma~\ref{l_oursharp}, but without the complications arising from non-square lattices;
we omit the details.
Since the event $H(R)$ depends only on variables $v_x$ for $x$ within a fixed distance
of $R$, taking $n$ large enough we may use Lemma~\ref{l_1dep} to deduce that $\thpw>0$:
the argument is
exactly that given in the last paragraph of
the proof of Theorem \ref{th_sitesq} in Subsection~\ref{ssite}; see also \cite{ourKesten}.
Similarly, we may use Lemma~\ref{l_kdepneg}
to deduce exponential decay for $p<1/2$, as in Section~\ref{bdecay}.
\end{proof}

\subsection{Random discrete Voronoi percolation}

Our final example is a discrete approximation of random Voronoi percolation.
Random Voronoi percolation, described below,
was introduced in the context of first-passage percolation by Vahidi-Asl and Wierman~\cite{VW}.
The critical probability, $1/2$, was established in~\cite{Voronoi}. The proof there is rather
long; the majority of the difficulties arise in attempting to compare Voronoi percolation
with a suitable discrete model, to which the method of \cite{ourKesten} can be applied.
Here we shall give a much simpler proof of a discrete result.

We start with $L=\Z^2$.
Given $0<\pi\le 1$, we select vertices of $L$ independently
at random, selecting each with probability $\pi$, to form a random set $L_\pi$.
Given $L_\pi$, we form the Voronoi cells associated to these points: for $z\in L_\pi$ let
\[
 V(z) = V_{L_\pi}(z) = \{x\in \RR^2: d(x,z)=\inf_{y\in L_\pi} d(x,y)\},
\]
where $d(.,.)$ is the Euclidean distance.
Thus $V(z)$ is the set of points in the plane at least as close to $z$ as to any other point $y$
of $L_\pi$. We include the boundary, obtaining with probability $1$
a set of closed convex polygons $V(z)$, $z\in L_\pi$,
that tile $\RR^2$.
We say that two cells $V(z_1)$, $V(z_2)$ are {\em weakly adjacent} if they share at least one point,
and {\em strongly adjacent} if they share an edge. These definitions may differ; indeed, they will
do so wherever four or more cells meet at a vertex. Given $L_\pi$ and $0<p<1$, we assign each Voronoi
cell a {\em state}, {\em open} or {\em closed},
taking cells to be open with probability $p$, independently of each other.
We write $\pp{p}$ for the associated probability measure.

A strong (weak) path of open cells is a sequence of open cells in
which each consecutive pair is strongly (weakly) adjacent.
The strong (weak) open component of the origin is the set of cells joined by a strong (weak) open path
to a cell containing the origin.

\begin{theorem}\label{th_dV}
Let $0<\pi\le 1$ and $0<p<1$ be given. If $p>1/2$, then with positive $\pp{p}$-probability the weak
component of the origin is infinite. If $p<1/2$, then there is a constant $a=a(\pi,p)>0$ such that the
$\pp{p}$-probability that the strong component of the origin contains
more than $n$ cells is at most $\exp(-an)$.
\end{theorem}

As $\pi\to 0$, after suitable rescaling $L_\pi$ converges to a Poisson process on $\RR^2$, and the Voronoi
tiling associated to $L_\pi$
approaches that associated to a Poisson process. In such a tiling, cells meet only three at a vertex,
so strong and weak connections coincide. Thus, for small $\pi$, the set-up considered in Theorem~\ref{th_dV}
is a good approximation to random Voronoi percolation, and the result strongly suggests that the critical
probability for random Voronoi percolation in the plane is $1/2$,
as shown in~\cite{Voronoi}. However, one cannot
just deduce this result (this would amount to an unjustified exchange of the order of two limits);
in fact, dealing with random Voronoi percolation is much harder.

\begin{proof}[Outline proof of Theorem~\ref{th_dV}]
Let us associate a random variable $v_z$ to each vertex $z$ of $L=\Z^2$.
We take $v_z=0$ if $z\notin L_\pi$, $v_z=+1$ if $z\in L_\pi$ and $V(z)$ is open,
and $v_z=-1$ if $z\in L_\pi$ and $V(z)$ is closed. Thus the $v_z$ are independent and identically
distributed, with $\Pr(v_z=i)=p_i$, where $p_{-1}=\pi(1-p)$, $p_0=1-\pi$ and $p_{+1}=\pi p$.
Let us say that a $x$ point of $\RR^2$ is {\em open} if it lies in an open cell.
Equivalently, $x$ is open if there is a $z\in L$ with $v_z=+1$ such that no $z'\in L$ 
with $d(x,z')<d(x,z)$ has $v_{z'}=-1$.
This event is increasing with respect to the $v_z$. Note that two cells $V(z)$, $V(z')$ are connected
by a weak open path if and only if there is a piecewise-linear path $P\subset \RR^2$ joining $z$
and $z'$ with every point of $P$ open.
Given a rectangle $R$, let us define horizontal and vertical open crossings of $R$ in terms
of such paths $P$.

We claim that the conditions of Theorem~\ref{th_wRSW} are satisfied.
Indeed, condition (i)
follows from our definition of openness for points of $\RR^2$. Condition (ii) is immediate -- our
set-up inherits the symmetries of the lattice $L=\Z^2$ we started from. (iii) is very easy to check:
for a fixed rectangle $R_1$, for large $\lambda$ it is very likely that every disc of radius
$\log\lambda$ centered within distance $\log\lambda$ of $\lambda R_1$ contains at least one point of
$L_\pi$; the expected number of such discs containing no points of $L_\pi$
tends to $0$ as $\lambda\to\infty$. It follows that with probability $1-o(1)$ the states of all points
in $\lambda R_1$ are determined by the variables $v_z$ for $z$ within distance $2\log\lambda=o(\lambda)$
of $\lambda R_1$; asymptotic independence follows immediately. For (iv),
very crudely, with probability $1-o(1)$ the length of a shortest
path crossing $\lambda R$
is at most $\lambda^3$, as all Voronoi cells meeting $\lambda R$ have
diameter at most $\log \lambda$, so there are at most $O(\lambda^2)$ such cells.

Let $R$ be any rectangle. Defining a point of $\RR^2$ to be {\em closed} if it lies in a closed cell (so 
some points are both open and closed, if they are in the boundary of an open cell and of a closed cell),
it is easy to check that either $R$ is crossed horizontally by an open path, or $R$ is crossed
vertically by a closed path, or both. (As usual, consider the topological boundary of the set of
open points in $R$ reachable by an open path from a point on the left-hand side of $R$.)
It follows that for any $\pi>0$ and any $p$ we have
$\pp{p}(H(R))+\pp{1-p}(V(R))\ge 1$. Thus, writing $R_{m,n}$ for an $m$ by $n$ rectangle,
$\pp{1/2}(H(R_{n,n}))\ge 1/2$ for all $n$.
Hence, by Theorem~\ref{th_wRSW}, there is a $c'>0$ such that
\begin{equation}\label{pp1}
 \pp{1/2}(H(R_{10n,n}))>c'
\end{equation}
for arbitrarily large $n$.

The rest of the argument is again similar to that in \cite{ourKesten} and in Section~\ref{bdecay}.
It suffices to show
that for any fixed $\pi$, $p>1/2$, $c''<1$ and $n_0$, there is an $n\ge n_0$ such that
\begin{equation}\label{pp2}
 \pp{p}(H(R_{6n,2n})) > c''.
\end{equation}
Then, recalling that open paths in $\RR^2$ correspond to weak paths of open cells,
the first statement of Theorem~\ref{th_dV} follows from Lemma~\ref{l_1dep} as usual
(see the proof of Theorem \ref{th_sitesq} in Subsection~\ref{ssite}, or~\cite{ourKesten}),
except that we must be a little careful defining the $1$-dependent measure: to achieve
$1$-dependence, we work with a modified form $G'(R)$ of the event $G(R)$, where $G'(R)$
depends only on the variables $v_z$ for $z$ within distance $n/2$, say, of the $6n$ by $2n$ rectangle $R$.
For $n_0$ large enough, we can find such a $G'(R)$ with probability close to that of $G(R)$;
the argument is as for asymptotic independence.  The same
technicality arises in the Voronoi setting; see \Sectionpfmain\ of~\cite{Voronoi}.
For the second statement, we use Lemma~\ref{l_kdepneg} as in Section~\ref{bdecay},
noting that if $C$ is a weak cycle of open cells, then no strong path of closed cells starts
inside and ends outside $C$.

To deduce \eqref{pp2} from \eqref{pp1}, we argue as in the proof of Lemma~\ref{l_oursharp}.
In this argument we have to overcome two additional
minor complications. Firstly, it is convenient to work in the product of three element
probability spaces, as above, so we need a version of Lemma~\ref{l_sharp'} that applies 
in this setting. Such a result is given in~\cite{Voronoi}; the proof is a very simple
modification of the proof of Theorem 3.2 of Friedgut and Kalai~\cite{FK}. Secondly, as the event $H(R)$
depends on points outside $R$, it is no longer quite true that the crossing probability of a rectangle in
$\RR^2$ and of the corresponding rectangle in the torus coincide. However, the difference
tends to zero as we enlarge the rectangle and torus in a constant ratio; the argument is the
same as for asymptotic independence above.
\end{proof}

\end{document}

%% file: sitedual.pstex_t
\begin{picture}(0,0)%
\includegraphics{sitedual.pstex}%
\end{picture}%
\setlength{\unitlength}{1302sp}%
\begingroup\makeatletter\ifx\SetFigFont\undefined%
\gdef\SetFigFont#1#2#3#4#5{%
  \reset@font\fontsize{#1}{#2pt}%
  \fontfamily{#3}\fontseries{#4}\fontshape{#5}%
  \selectfont}%
\fi\endgroup%
\begin{picture}(13605,10515)(211,-9877)
\put(226,239){\makebox(0,0)[b]{\smash{{\SetFigFont{10}{12.0}{\rmdefault}{\mddefault}{\updefault}{\color[rgb]{0,0,0}$x$}%
}}}}
\put(13651,-9736){\makebox(0,0)[b]{\smash{{\SetFigFont{10}{12.0}{\rmdefault}{\mddefault}{\updefault}{\color[rgb]{0,0,0}$z$}%
}}}}
\put(301,-9736){\makebox(0,0)[b]{\smash{{\SetFigFont{10}{12.0}{\rmdefault}{\mddefault}{\updefault}{\color[rgb]{0,0,0}$w$}%
}}}}
\put(13801,239){\makebox(0,0)[b]{\smash{{\SetFigFont{10}{12.0}{\rmdefault}{\mddefault}{\updefault}{\color[rgb]{0,0,0}$y$}%
}}}}
\end{picture}%

%% file: l_X3.pstex_t
\begin{picture}(0,0)%
\includegraphics{l_X3.pstex}%
\end{picture}%
\setlength{\unitlength}{2171sp}%
\begingroup\makeatletter\ifx\SetFigFont\undefined%
\gdef\SetFigFont#1#2#3#4#5{%
  \reset@font\fontsize{#1}{#2pt}%
  \fontfamily{#3}\fontseries{#4}\fontshape{#5}%
  \selectfont}%
\fi\endgroup%
\begin{picture}(6644,5444)(1779,-5783)
\put(2701,-4786){\makebox(0,0)[b]{\smash{{\SetFigFont{10}{12.0}{\rmdefault}{\mddefault}{\updefault}{\color[rgb]{0,0,0}$P_1$}%
}}}}
\put(1876,-1261){\makebox(0,0)[lb]{\smash{{\SetFigFont{10}{12.0}{\rmdefault}{\mddefault}{\updefault}{\color[rgb]{0,0,0}$R_3$}%
}}}}
\put(1876,-4861){\makebox(0,0)[lb]{\smash{{\SetFigFont{10}{12.0}{\rmdefault}{\mddefault}{\updefault}{\color[rgb]{0,0,0}$R_1$}%
}}}}
\put(1876,-3061){\makebox(0,0)[lb]{\smash{{\SetFigFont{10}{12.0}{\rmdefault}{\mddefault}{\updefault}{\color[rgb]{0,0,0}$R_2$}%
}}}}
\put(5926,-2986){\makebox(0,0)[lb]{\smash{{\SetFigFont{10}{12.0}{\rmdefault}{\mddefault}{\updefault}{\color[rgb]{0,0,0}$P'$}%
}}}}
\put(2626,-3211){\makebox(0,0)[b]{\smash{{\SetFigFont{10}{12.0}{\rmdefault}{\mddefault}{\updefault}{\color[rgb]{0,0,0}$P_2$}%
}}}}
\put(2626,-1186){\makebox(0,0)[b]{\smash{{\SetFigFont{10}{12.0}{\rmdefault}{\mddefault}{\updefault}{\color[rgb]{0,0,0}$P_3$}%
}}}}
\put(2851,-3511){\makebox(0,0)[lb]{\smash{{\SetFigFont{10}{12.0}{\rmdefault}{\mddefault}{\updefault}{\color[rgb]{0,0,0}$v$}%
}}}}
\put(5551,-661){\makebox(0,0)[lb]{\smash{{\SetFigFont{10}{12.0}{\rmdefault}{\mddefault}{\updefault}{\color[rgb]{0,0,0}$R$}%
}}}}
\end{picture}%

%% file: RSRk.pstex_t
\begin{picture}(0,0)%
\includegraphics{RSRk.pstex}%
\end{picture}%
\setlength{\unitlength}{2171sp}%
\begingroup\makeatletter\ifx\SetFigFont\undefined%
\gdef\SetFigFont#1#2#3#4#5{%
  \reset@font\fontsize{#1}{#2pt}%
  \fontfamily{#3}\fontseries{#4}\fontshape{#5}%
  \selectfont}%
\fi\endgroup%
\begin{picture}(5444,3624)(279,-3673)
\put(3001,-1936){\makebox(0,0)[b]{\smash{{\SetFigFont{10}{12.0}{\rmdefault}{\mddefault}{\updefault}{\color[rgb]{0,0,0}$S$}%
}}}}
\put(4426,-1711){\makebox(0,0)[b]{\smash{{\SetFigFont{10}{12.0}{\rmdefault}{\mddefault}{\updefault}{\color[rgb]{0,0,0}$R$}%
}}}}
\put(1351,-1711){\makebox(0,0)[b]{\smash{{\SetFigFont{10}{12.0}{\rmdefault}{\mddefault}{\updefault}{\color[rgb]{0,0,0}$R'$}%
}}}}
\end{picture}%